\documentclass[11pt]{article}
\usepackage{amsmath,fullpage,amssymb,amsthm,hyperref,enumerate,shuffle}
\usepackage{tikz,color}

\newtheorem{theorem}{Theorem}

\theoremstyle{remark}
\newtheorem{remark}[theorem]{Remark}  

\theoremstyle{definition}
\newtheorem{example}{Example}  

\newcommand{\D}{\mathcal{D}}
\newcommand{\M}{\mathcal{M}}
\newcommand{\F}{\overline{\mathcal{D}}}
\newcommand{\f}{\overline{d}}
\renewcommand{\S}{\mathcal{S}}
\newcommand{\I}{\mathbf{I}}   
\newcommand{\ul}{\underline}
\newcommand{\bij}{\phi}
\newcommand{\bijf}{\varphi}

\newcommand{\minus}{\!\setminus\!}

\date{}

\title{Bijections for restricted inversion sequences and permutations with fixed points}

\author{Sergi Elizalde\thanks{Department of Mathematics, Dartmouth College, Hanover, NH 03755. \texttt{sergi.elizalde@dartmouth.edu}}}

\begin{document}

\maketitle

\begin{abstract}
We provide a bijective proof of a formula of Auli and the author expressing the number of inversion sequences with no three consecutive equal entries in terms of the number of non-derangements, that is, permutations with fixed points. Additionally, we give bijective proofs of two simple recurrences for the number of non-derangements.
\end{abstract}

\section{Introduction}

Let $\S_n$ denote the set of permutations of $[n]=\{1,2,\dots,n\}$.  
A {\em fixed point} of $\pi\in\S_n$ if an element $i\in[n]$ such that $\pi(i)=i$. A {\em derangement} is a permutation without fixed points. Let $\D_n$ by the set of derangements in $\S_n$, and let $d_n=|\D_n|$. A permutation that is not a derangement, that is, has at least one fixed point, will be called a {\em non-derangement}. Denote the set of non-derangements in $\S_n$ by $\F_n$, and let $\f_n=|\F_n|=n!-d_n$. 

An {\em inversion sequence} is an integer sequence $e_{1}e_{2}\dots e_{n}$ such that $0\leq e_{i}<i$ for each $i$. 
Let $\I_n$ denote the set of inversion of length $n$. 
Inversion sequences provide a useful way to encode permutations. Several such encodings ---namely, bijections between $\S_n$ and $\I_n$--- have been extensively studied, sometimes under the name of {\em inversion tables}~\cite[Prop.\ 1.3.12]{Stanley} or {\em Lehmer codes}. Each entry $e_i$ counts certain inversions (that is, pairs of elements that appear out of order in the permutation), so that the sum of the entries of the inversion sequence equals the number of inversions of the permutation that it encodes.

Motivated by the vast literature on pattern avoidance in permutations, there has been growing interest in the study pattern avoidance in inversion sequences, including classical patterns~\cite{MartinezSavageI,MansourShattuck,MartinezSavageII}, consecutive patterns~\cite{AuliElizalde,AuliElizaldeII} and vincular patterns~\cite{LinYan,AuliElizaldeIII}. In~\cite{AuliElizalde}, Auli and the author enumerate inversion sequences avoiding each consecutive pattern of length~3, and in particular, avoiding the pattern $\ul{000}$.
An inversion sequence is said to avoid $\ul{000}$ if it does not contain three consecutive equal entries, that is, there does not exist $i\in[n-2]$ such that $e_i=e_{i+1}=e_{i+2}$. Let $$\I_n(\ul{000})=\{e\in\I_n:e\text{ avoids }\ul{000}\}.$$
For example, $\I_4(\ul{000})$ consists of the 19 inversion sequences $0010,0011,0012,0013,0020,0021,0022$, $0023,0100, 0101,0102,0103,0110,
0112,0113,0120,0121,0122,0123$.

In~\cite[Cor.\ 3.3]{AuliElizalde}, the authors give the following formula expressing the number of inversion sequences that avoid $\ul{000}$ in terms of the number of derangements. 

\begin{theorem}[\cite{AuliElizalde}]\label{thm:aulielizalde}
For $n\geq 1$,
\begin{equation}\label{eq:aulielizalde}
|\I_{n}(\underline{000})|=\frac{(n+1)!-d_{n+1}}{n}.
\end{equation}
\end{theorem}

This formula corresponds to sequence A052169 in the Online Encyclopedia of Integer Sequences~\cite{OEIS}, and its first few terms are $1, 2, 5, 19, 91, 531, 3641,\dots$. The proof in~\cite{AuliElizalde} is by induction on $n$, and it relies on the recurrences
$$|\I_{n}(\underline{000})|=(n-1)|\I_{n-1}(\underline{000})|+(n-2)|\I_{n-2}(\underline{000})|$$
and
\begin{equation}
\label{eq:dn} d_{n}=(n-1)(d_{n-1}+d_{n-2})
\end{equation}
for $n\ge2$. Finding a bijective proof of Theorem~\ref{thm:aulielizalde} is left as a open problem in~\cite{AuliElizalde}.
We will provide such proof in Section~\ref{sec:bij}.

In Section~\ref{sec:F} we introduce two natural recurrences for the non-derangement numbers $\f_n$, which are analogous to well-known recurrences for the derangement numbers, and we provide bijective proofs of them.

\section{A bijective proof of Theorem~\ref{thm:aulielizalde}}\label{sec:bij}

As noted by Peter Winkler (personal communication, May 5, 2020), the right-hand side of Equation~\eqref{eq:aulielizalde} can be interpreted as the cardinality of the set $\F_n\sqcup\F_{n-1}$, where $\sqcup$ denotes a disjoint union. Indeed, using Equation~\eqref{eq:dn},
\begin{align*}
\frac{(n+1)!-d_{n+1}}{n}&=\frac{(n+1)!-n(d_{n}+d_{n-1})}{n}=(n+1)(n-1)!-d_n-d_{n-1}\\
&=n!-d_n+(n-1)!-d_{n-1}
=\f_n+\f_{n-1}.
\end{align*}
Our proof of Theorem~\ref{thm:aulielizalde} will consist of a bijection
$$\bij:\I_n(\ul{000})\to \F_n\sqcup\F_{n-1}.$$

\subsection{The bijection}

In order to describe the map $\bij$, first we introduce some notation. If $1\le a,b\le n$, denote by $(a,b)$ the permutation in $\S_n$ that switches $a$ and $b$; this permutation is a transposition if $a\neq b$, and it is the identity permutation if $a=b$. In particular, if $\sigma=\sigma(1)\sigma(2)\dots\sigma(n)\in\S_n$, then the product $(a,b)\sigma$ denotes the permutation obtained from the one-line notation of $\sigma$ by switching the entries $a$ and $b$ if they are different. For, example, $(3,5)32514=52314$ and $(2,2)32514=32514$.
We will also use the fact that $\S_{n-1}$ can be viewed as the subset of $\S_n$ consisting of those permutations where $n$ is a fixed point.

Let $n\ge1$, and let $e\in\I_n(\ul{000})$.
The first step in the construction of $\phi(e)$ is to encode $e$ as a word $w=w_2\dots w_n$ as follows. For $2\le k\le n$, let
$$w_k=\begin{cases} R & \text{if } e_k=e_{k-1}, \\
e_k & \text{if } e_k>e_{k-1},\\
e_k+1 & \text{if } e_k<e_{k-1}. \end{cases}$$
This encoding is a bijection between $\I_n(\ul{000})$ and the set of words $w=w_2\dots w_n$ with $w_k\in [k-1]\cup\{R\}$ (for $2\le k\le n$) not containing two consecutive $R$s.

Next we read $w$ from left to right and build a sequence of permutations $\sigma_1,
\sigma_2,\dots, \sigma_n$, where $\sigma_k\in \F_k\sqcup\F_{k-1}$ for all $k$;
 specifically, $\sigma_k\in\F_k$ if $w_k\neq R$, and $\sigma_k\in\F_{k-1}$ if $w_k=R$.
Set $\sigma_1=1\in\F_1$. For each $k$ from $2$ to $n$, repeat the following step. 
If $w_k=R$, let $\sigma_k=\sigma_{k-1}\in\F_{k-1}$. Otherwise, let
$$\sigma_k=\begin{cases}
(w_k,k)\sigma_{k-1} & \text{if $w_{k-1}\neq R$ and $\sigma_{k-1}\in\F_{k-1}$ has fixed points other than $w_k$},\\
(w_k,k-1)\sigma_{k-1} & \text{otherwise}, 
\end{cases}
$$
where we view $\sigma_{k-1}$ as an element of $\F_k$, and thus $\sigma_k\in\F_k$.

Finally, we define $\bij(e)=\sigma_n$.

\begin{example}
If $e=001322$, then $w=R133R$. The computation of $\sigma_k$ for $k$ from $1$ to $6$ is shown in Table~\ref{tab:ex1}, giving
$\bij(e)=\sigma_6=21543\in\F_6\sqcup\F_5$.

If $e=0102230$, then $w=112R31$. The computations in Table~\ref{tab:ex2} give
$\bij(e)=\sigma_7=2574361\in\F_7\sqcup\F_6$.
\end{example}

\begin{table}[h]
$$\begin{array}{c|c|c|c}
k & w_k & \sigma_k \text{ in one-line notation} & \sigma_k \text{ in cycle notation} \\ \hline
1 & & 1 & (1)\\
2 & R & 1 & (1)\\
3 & 1 & (12)123=213 & (2,1)(3)\\
4 & 3 & (3,3)2134=2134 & (2,1)(3)(4)\\
5 & 3 & (3,5)21345=21543 & (2,1)(5,3)(4)\\
6 & R & 21543 & (2,1)(5,3)(4) 
\end{array}
$$
\caption{The computation of $\bij(001322)=21543$.}
\label{tab:ex1}
\end{table}

\begin{table}[h]
$$\begin{array}{c|c|c|c}
k & w_k &  \sigma_k \text{ in one-line notation} & \sigma_k \text{ in cycle notation} \\ \hline
1 & & 1 & (1)\\ 
2 & 1 & (1,1)12=12 & (1)(2)\\
3 & 1 & (1,3)123=321 & (3,1)(2) \\
4 & 2 & (2,3)3214=2314 & (2,3,1)(4) \\
5 & R & 2314 & (2,3,1)(4) \\
6 & 3 & (3,5)231456=251436 & (2,5,3,1)(4)(6)\\
7 & 1 & (1,7)2514367=2574361 & (2,5,3,7,1)(4)(6)
\end{array}
$$
\caption{The computation of $\bij(0102230)=2574361$.}
\label{tab:ex2}
\end{table}

\subsection{The map $\bij$ in cycle notation}
It is sometimes convenient to describe the construction of $\sigma_1,\sigma_2,\dots, \sigma_n$ in cycle notation, where we write permutations as products of disjoint cycles. 
In this case, $\sigma_1=(1)\in\F_1$, and for each $k$ from $2$ to $n$, we repeat the following step. If $w_k=R$, let $\sigma_k=\sigma_{k-1}\in\F_{k-1}$. Otherwise,  $\sigma_k\in\F_k$ is obtained from the cycle notation of $\sigma_{k-1}$ as follows:
\begin{itemize}
\item if $w_{k-1}=R$, insert $k-1$ right before $w_k$ in the same cycle, and add a new fixed point $(k)$;
\item if $w_{k-1}\neq R$ and $\sigma_{k-1}\in\F_{k-1}$ has fixed points other than $w_k$, insert $k$ right before $w_k$ in the same cycle,
\item otherwise (that is, if $w_{k-1}\neq R$ and $w_k$ is the only fixed point of $\sigma_{k-1}\in\F_{k-1}$), add a new fixed point $(k)$, and if $w_k\neq k-1$, remove $(w_k)$ and insert $w_k$ right before $k-1$ in the same cycle.
\end{itemize}
The right column of Tables~\ref{tab:ex1} and~\ref{tab:ex2} shows examples of this construction in cycle notation.

\subsection{The inverse map}
To show that $\bij$ is indeed a bijection, let us describe its inverse $\bij^{-1}:\F_n\sqcup\F_{n-1}\to\I_n(\ul{000})$.
Given $\pi\in\F_n\sqcup\F_{n-1}$, set $\sigma_n=\pi$. We will describe permutations
$\sigma_{n-1}, \sigma_{n-2},\dots, \sigma_1$ in cycle notation, while building a word $w$ from right to left.
For $k$ from $n$ to $2$, repeat the following step. If $\sigma_k\in\F_{k-1}$, let $w_k=R$ and $\sigma_{k-1}=\sigma_k$. Otherwise (that is, if $\sigma_k\in\F_{k}$) proceed as follows.
\begin{itemize}
\item If $k$ is not a fixed point of $\sigma_k$, let $w_k=\sigma_k(k)$, and let $\sigma_{k-1}$ be the permutation obtained by removing $k$ from the cycle notation of $\sigma_k$.
\item Otherwise, remove $(k)$ from the cycle notation of $\sigma_k$, and then:
\begin{itemize}
\item if removing $k-1$ from the cycle notation leaves any fixed points, let $\sigma_{k-1}\in\F_{k-2}$ be the resulting permutation, and let $w_k=\sigma_k(k-1)$;
\item otherwise, let $w_k=\sigma_k^{-1}(k-1)$
, and let 
$\sigma_{k-1}\in\F_{k-1}$ be the permutation obtained by removing $w_k$ from its current cycle and creating a fixed point $(w_k)$. (Note that this produces no change if $w_k=k-1$, since in this case $w_k$ is already a fixed point.) 
\end{itemize}
\end{itemize}

From $w$, we obtain the inversion sequence $e=\bij^{-1}(\pi)$ by letting $e_1=0$ and, for $2\le k\le n$, letting
$$e_k=\begin{cases} e_{k-1} & \text{if } w_k=R, \\
w_k & \text{if } w_k>e_{k-1},\\
w_k-1 & \text{if } w_k\le e_{k-1}. \end{cases}$$

\subsection{A related bijection}

It is suggested in~\cite{AuliElizalde} that, to prove Theorem~\ref{thm:aulielizalde} bijectively, one could give a bijection between the sets $[n]\times \I_{n}(\underline{000})$ and $\F_{n+1}$, which have cardinalities $n |\I_n(\ul{000})|$ and $\f_{n+1}=(n+1)!-d_{n+1}$, respectively. 
Let us show how to use $\bij$ to provide such a bijection.
Given $a\in[n]$ and $e=e_1e_2\dots e_n\in\I_{n}(\underline{000})$, define $e_{n+1}=a$ if $a>e_{n}$, and $e_{n+1}=a-1$ if $a\le e_n$.
Then $e'=e_1e_2\dots e_ne_{n+1}\in\I_{n+1}(\underline{000})$, and $\bij(e')\in\F_{n+1}$. The map $(a,e)\mapsto \bij(e')$ gives the desired bijection. Indeed, our construction provides a composition of bijections
$$\begin{array}{ccccc}
[n]\times \I_{n}(\underline{000}) & \to & \{e_1e_2\dots e_ne_{n+1}\in\I_{n+1}(\underline{000}):e_n\neq e_{n+1}\} & \to & \F_{n+1}\\
(a,e) & \mapsto & e' & \mapsto & \bij(e').
\end{array}$$

An example of this composition for $n=6$ is $(1,010223)\mapsto 0102230 \mapsto 2574361$.

\section{Bijective proofs of non-derangement recurrences}\label{sec:F}

\subsection{First recurrence}

Equation~\eqref{eq:dn} is a well-known derangement recurrence (see e.g.\ \cite[Eq.~(2.14)]{Stanley}) having a simple combinatorial proof. Indeed, letting $\D_n$ be the set of derangements in $\S_n$, a bijection $$\D_n\to[n-1]\times\D_{n-1}\sqcup[n-1]\times\D_{n-2}$$ is obtained by mapping $\pi\in\D_n$ to the pair $(\pi(n),\pi')$, where $\pi'$ is obtained from $\pi$ as follows: if $n$ belongs to a $2$-cycle in $\pi$, remove this $2$-cycle; otherwise, remove $n$ from the cycle notation of $\pi$.

Equation~\eqref{eq:dn} 
implies that non-derangement numbers also satisfy the recurrence
\begin{equation}
\label{eq:fn1} \f_{n}=(n-1)(\f_{n-1}+\f_{n-2})
\end{equation}
for $n\ge2$, with initial terms $\f_0=0$, $\f_1=1$ (compare to $d_0=1$, $d_1=0$).
Indeed, 
\begin{align*}\f_n&=n!-d_n=n!-(n-1)(d_{n-1}+d_{n-2})=n!-(n-1)\left((n-1)!-\f_{n-1}+(n-2)!-\f_{n-2}\right)\\
&=n!-(n-1)\left((n-1)!+(n-2)!\right)-\f_{n-1}-\f_{n-2}=(n-1)(\f_{n-1}+\f_{n-2}).\end{align*}
The values of the sequence $\f_n$ for $1\le n\le 7$ are $1, 1, 4, 15, 76, 455, 3186$. This is sequence A002467 in~\cite{OEIS}.

Next we provide a direct combinatorial proof of Equation~\eqref{eq:fn1}, by describing a bijection 
$$\bijf:\F_n\to[n-1]\times\F_{n-1}\sqcup[n-1]\times\F_{n-2}.$$
This bijection is implicitly used in the construction of $\bij$ in Section~\ref{sec:bij}.
For $\pi\in\S_n$ and $S\subseteq[n]$, denote by $\pi\minus S$ the permutation of $[n]\setminus S$ obtained from the cycle notation of $\pi$ by removing the elements in $S$. For example, if $\pi=(1,6)(2,5,3)(4)$, then $\pi\minus\{6,4\}=(1)(2,5,3)$. 

For $\pi\in\F_n$, define
\begin{equation}\label{def:bijf}
\bijf(\pi)=\begin{cases}
(\pi(n),\pi\minus\{n\}) & \text{if }\pi(n)\neq n,\\
(\pi(n-1),\pi\minus\{n,n-1\}) & \text{if }\pi(n)=n\text{ and }\pi\minus\{n,n-1\}\in\F_{n-2},\\
(\pi^{-1}(n-1),(\pi^{-1}(n-1))\,\pi\minus\{n,\pi^{-1}(n-1)\}) & \text{otherwise.}
\end{cases}
\end{equation}
An example is given in Table~\ref{tab:bijf}, which lists $\bij(\pi)$ for each $\pi\in\F_4$.

\begin{table}[htb]
$$\begin{array}{c|c}
\pi\in\F_4 & \bijf(\pi) \\ \hline
(1)(2)(3)(4) & (3,(1)(2)) \\
(1,2)(3)(4) & (3,(3)(1,2)) \\
(1,3)(2)(4) & (1,(1)(2)) \\
(1,4)(2)(3) & (1,(1)(2)(3)) \\
(1)(2,3)(4) & (2,(1)(2)) \\
(1)(2,4)(3) & (2,(1)(2)(3)) \\
(1)(2)(3,4) & (3,(1)(2)(3)) \\
(1,2,3)(4) & (2,(2)(1,3)) \\
(1,3,2)(4) & (1,(1)(3,2)) \\
(1,2,4)(3) & (1,(1,2)(3)) \\
(1,4,2)(3) & (2,(1,2)(3)) \\
(1,3,4)(2) & (1,(1,3)(2)) \\
(1,4,3)(2) & (3,(1,3)(2)) \\
(1)(2,3,4) & (2,(1),(2,3)) \\
(1)(2,4,3) & (3,(1),(2,3))
\end{array}$$
\caption{The bijection $\bijf:\F_4\to [3]\times\F_{3}\sqcup[3]\times\F_{2}$.}
\label{tab:bijf}
\end{table}

The inverse map $\bijf^{-1}:[n-1]\times\F_{n-1}\sqcup[n-1]\times\F_{n-2}\to\F_n$ can be described as follows. 
For $(i,\sigma)\in[n-1]\times\F_{n-1}\sqcup[n-1]\times\F_{n-2}$, let 
$$\bijf^{-1}((i,\sigma))=\begin{cases}
(i,n)\sigma & \text{if $\sigma\in\F_{n-1}$ and $\sigma$ has fixed points other than $i$},\\
(i,n-1)\sigma & \text{otherwise}, 
\end{cases}
$$
where, in the expressions $(a,b)\sigma$ on the right-hand side, we view $\sigma$ as a permutation in $\S_n$.

\begin{remark}
An alternative bijection is obtained by replacing the third case in Equation~\eqref{def:bijf} with 
$$(\pi(n-1),(\pi(n-1))\,\pi\minus\{n,\pi(n-1)\}).$$
This variation still describes a bijection from $\F_n$ to $[n-1]\times\F_{n-1}\sqcup[n-1]\times\F_{n-2}$, but the description of its inverse requires an additional separate case for pairs $(i,\sigma)$ where $\sigma\in\F_{n-1}$ and $i$ is the only a fixed point of $\sigma$.
\end{remark}

\subsection{Second recurrence}
Another well-known recurrence (see e.g.\ \cite[Eq.~(2.13)]{Stanley}) for the derangement numbers is
\begin{equation}
\label{eq:dn2}
d_n=nd_{n-1}+(-1)^n
\end{equation}
for $n\ge1$. Combinatorial proofs of this recurrence, which require considerably more work than for Equation~\eqref{eq:dn}, have appeared in~\cite{Desarmenien,Remmel,Wilf,Rak,Benjamin,Elizalde_der}. From the fact that $\f_n=n!-d_n$, it follows immediately from Equation~\eqref{eq:dn2} that
\begin{equation}\label{eq:fn2}
\f_n=n\f_{n-1}-(-1)^n
\end{equation}
for $n\ge1$. Next we provide a direct bijective proof of Equation~\eqref{eq:fn2}.

Let us start by giving a combinatorial interpretation of $n\f_{n-1}$.  This is the number of permutations in $\S_n$ with a {\em marked} (i.e., distinguished) fixed point, and at least one unmarked fixed point. Indeed, there are $n$ ways to choose an element of $[n]$ to be the marked fixed point, and $d_{n-1}$ to choose a derangement of the remaining $n-1$ elements. Let $\M_n$ be the set of such marked permutations. We will write them in cycle notation with the marked fixed point underlined, such as in $(1,5)(\ul{2})(3,7,6)(4)\in\M_7$.

To prove Equation~\eqref{eq:fn2} combinatorially, we describe a bijection 
$\theta:\F^\ast_n\to \M^\ast_n$, where $\F^\ast_n=\F_n\setminus\{(1,2)(3,4)\dots(n-2,n-1)(n)\}$ and $\M^\ast_n=\M_n$ when $n$ is odd, and 
$\F^\ast_n=\F_n$ and $\M^\ast_n=\M_n\setminus\{(\ul{1})(2,3)\dots(n-2,n-1)(n)\}$ when $n$ is even.

Given $\pi\in\F^\ast_n$ written in cycle notation, consider three cases:
\begin{enumerate}[A.]
\item If $(n)$ is a fixed point of $\pi$ but not the only one, let $\theta(\pi)$ be obtained from $\pi$ by marking the fixed point $(\ul{n})$.
\item If $(n)$ is the only fixed point of $\pi$, write the cycles of $\pi$ beginning with their smallest element, and listed by increasing first element. Let $k$ be the largest non-negative integer such that $\pi$ starts with $(1,2)(3,4)\dots(2k-1,2k)$, and note that $0\le k<\frac{n-1}{2}$.
Let $\theta(\pi)$ be obtained as follows:
\begin{enumerate}[i.]
\item If the cycle of $\pi$ containing $2k+1$ has at least 3 elements, change the first $k+1$ cycles of $\pi$ as follows:
\begin{align*}
\pi&=(1,2)(3,4)\dots(2k-1,2k)(2k+1,a_1,a_2,\dots,a_j)\dots(n)\\
\theta(\pi)&=(\ul{1})(2,3)(4,5)\quad\dots\quad(2k,a_1)(2k+1,a_2,\dots,a_j)\dots(n)
\end{align*}
Note that in the case $k=0$, this construction gives $\theta(\pi)=(\ul{a_1})(1,a_2,\dots,a_j)\dots(n)$.
\item If the cycle of $\pi$ containing $2k+1$ has 2 elements, change the first $k+2$ cycles of $\pi$ as follows:
\begin{align*}
\pi&=(1,2)(3,4)\dots(2k-1,2k)(2k+1,a_1)(2k+2,a_2,\dots,a_j)\dots(n) \\
\theta(\pi)&=(\ul{1})(2,3)(4,5)\,\quad\dots\quad\,(2k,2k+1)(2k+2,a_1,a_2,\dots,a_j)\dots(n)
\end{align*}
\end{enumerate}
\item If $\pi(n)\neq n$, let $\theta(\pi)$ be obtained from $\pi$ by removing $\pi(n)$ from its cycle and creating a new marked fixed point $(\ul{\pi(n)})$.
\end{enumerate}

Table~\ref{tab:theta} lists $\theta(\pi)$ for each $\pi\in\F^\ast_4=\F_4$.

\begin{table}[htb]
$$\begin{array}{c|c}
\pi\in\F^\ast_4 & \theta(\pi)\in\M^\ast_4 \\ \hline
(1)(2)(3)(4) &  (1)(2)(3)(\ul{4})\\
(1,2)(3)(4) &  (1,2)(3)(\ul{4}) \\
(1,3)(2)(4) & (1,3)(2)(\ul{4}) \\
(1,4)(2)(3) &  (\ul{1})(4)(2)(3) \\
(1)(2,3)(4) & (1)(2,3)(\ul{4})  \\
(1)(2,4)(3) &  (1)(\ul{2})(3)(4) \\
(1)(2)(3,4) &  (1)(2)(\ul{3})(4)\\
(1,2,3)(4) &  (1,3)(\ul{2})(4) \\
(1,3,2)(4) & (1,2)(\ul{3})(4) \\
(1,2,4)(3) &  (\ul{1})(2,4)(3) \\
(1,4,2)(3) & (1,4)(\ul{2})(3) \\
(1,3,4)(2) &  (\ul{1})(2)(3,4) \\
(1,4,3)(2) &  (1,4)(2)(\ul{3}) \\
(1)(2,3,4) &  (1)(\ul{2})(3,4) \\
(1)(2,4,3) &  (1)(2,4)(\ul{3})
\end{array}$$
\caption{The bijection $\theta:\F_4\to \M_4\setminus\{(\ul{1})(2,3)(4)\}$.}
\label{tab:theta}
\end{table}

In case B above, which corresponds to permutations consisting of a derangement $\pi'\in\D_{n-1}$ together with the fixed point $(n)$, the description of $\theta(\pi)$ is equivalent to applying the bijection $\psi$ from~\cite{Elizalde_der} (between derangements and permutations with one fixed point) to $\pi'\in\D_{n-1}$, and then marking the fixed point of $\psi(\pi)$. Below are some more examples of this case.

\begin{example}
If $\pi=(1,2)(3,5,6,4)(7)$, we are in case B.i with $k=1$, and $\theta(\pi)=(\ul{1})(2,5)(3,6,4)(7)$. 
If $\pi=(1,2,4)(3,5)(6)$, we are in case B.i with $k=0$, and $\theta(\pi)=(\ul{2})(1,4)(3,5)(7)$. 
If $\pi=(1,2)(3,6)(4,5)(7)$, we are in case B.ii with $k=1$, and $\theta(\pi)=(\ul{1})(2,3)(4,6,5)(7)$.
\end{example}

To show that $\theta$ is indeed a bijection, next we describe its inverse map. Given $\sigma\in\M_n$, we can recover $\theta^{-1}(\sigma)$ in each case as follows.

\begin{enumerate}[A.]
\item If $(\ul{n})$ is the marked fixed point of $\sigma$, let $\theta^{-1}(\sigma)$ be obtained by unmarking it.
\item If $(n)$ is the only unmarked fixed point of $\sigma$, suppose $(\ul{\ell})$ is the marked fixed point. If $\ell\neq 1$, let $\theta^{-1}(\sigma)$ be obtained from $\sigma$ by removing  $(\ul{\ell})$ and inserting $\ell$ right after $1$ in the same cycle.
If $\ell=1$, write the cycles of $\sigma$ beginning with their smallest element, and listed by increasing first element. 
Let $k'$ be the largest positive integer such that $\sigma$ starts with $(\ul{1})(2,3)(4,5)\dots(2k'-2,2k'-1)$, and note that $1\le k'<\frac{n-1}{2}$.
\begin{enumerate}[i.]
\item If the cycle of $\sigma$ containing $2k'$ has 2 elements, change the first $k'+2$ cycles of $\sigma$ as follows:
\begin{align*}
\sigma&=(\ul{1})(2,3)(4,5)\dots(2k'-2,2k'-1)(2k',a_1)(2k'+1,a_2,\dots,a_j)\dots(n)\\
\theta^{-1}(\sigma)&=(1,2)(3,4)\, \ \qquad\dots\qquad\ \, (2k'-1,2k')(2k'+1,a_1,a_2,\dots,a_j)\dots(n)
\end{align*}
\item If the cycle of $\sigma$ containing $2k'$ has at least 3 elements, change the first $k'+1$ cycles of $\sigma$ as follows:
\begin{align*}
\sigma&=(\ul{1})(2,3)(4,5)\,\quad\dots\quad\,(2k'-2,2k'-1)(2k',a_1,a_2,\dots,a_j)\dots(n)\\
\pi&=(1,2)(3,4)\dots(2k'-3,2k'-2)(2k'-1,a_1)(2k',a_2,\dots,a_j)\dots(n)
\end{align*}
\end{enumerate}

\item Otherwise
, suppose $(\ul{\ell})$ is the marked fixed point, and let $\theta^{-1}(\sigma)$ be obtained by removing $(\ul{\ell})$ from $\sigma$ and inserting $\ell$ right after $n$ in the same cycle of $\sigma$.
\end{enumerate}

\bibliographystyle{plain}
\bibliography{bijectionderangements}

\end{document}